\documentclass[11pt]{amsart}
\usepackage{amssymb}
\newtheorem{theorem}{Theorem}[section]
\newtheorem{lemma}[theorem]{Lemma}
\newtheorem{proposition}[theorem]{Proposition}
\newtheorem{corollary}[theorem]{Corollary}
\newtheorem{definition}[theorem]{Definition}

\newtheorem{question}[theorem]{Question}
\newcommand{\C}{\mathbb{C}}
\newcommand{\Z}{\mathbb{Z}}
\newcommand{\R}{\mathbb{R}}
\newcommand{\CP}{\mathbb{CP}}
\title[Symplectic 4-manifolds, plane curves and braid groups]%
{Some open questions about symplectic 4-manifolds,
singular plane curves, and\\ braid group factorizations}
\author{Denis Auroux}
\thanks{Partially supported by NSF grant DMS-0244844.}
\address{Department of Mathematics, M.I.T., Cambridge MA 02139, USA}
\email{auroux@math.mit.edu}

\begin{document}
\begin{abstract}
The topology of symplectic 4-manifolds is related to that of singular
plane curves via the concept of branched covers. Thus, various
classification problems concerning symplectic 4-manifolds can be
reformulated as questions about singular plane curves.
Moreover, using braid monodromy, these can in turn be
reformulated in the language of braid group factorizations.
While the results mentioned in this paper are not new, we hope that they
will stimulate interest in these questions, which remain essentially wide
open.
\end{abstract}

\maketitle

\section{Introduction}

An important problem in 4-manifold topology is to understand which manifolds
carry symplectic structures (i.e., closed non-degenerate 2-forms), and to
develop invariants that can distinguish symplectic manifolds. Additionally,
one would like to understand to what extent the category of symplectic
manifolds is richer than that of K\"ahler (or complex projective) manifolds.
For example, one would like to identify a set of surgery operations that can
be used to turn an arbitrary symplectic 4-manifold into a K\"ahler manifold,
or two symplectic 4-manifolds with the same classical topological invariants
(fundamental group, Chern numbers, ...) into each other.

Similar questions may be asked about singular curves inside, e.g., the
complex projective plane. The two types of questions are related to each
other via symplectic branched covers.
A branched cover of a symplectic 4-manifold with a (possibly singular)
symplectic branch curve carries a natural symplectic structure. Conversely,
every compact symplectic 4-manifold is a branched cover of $\CP^2$,
with a branch curve presenting nodes (of both orientations) and
complex cusps as its only singularities. 

In the language of branch curves, the failure of most symplectic manifolds
to admit integrable complex structures translates into the failure of most
symplectic branch curves to be isotopic to complex curves. While the
symplectic isotopy problem has a negative answer for plane curves with
cusp and node singularities, it is interesting to investigate this failure
more precisely. Various partial results have been obtained recently about
situations where isotopy holds (for smooth curves; for curves of low
degree), and about isotopy up to stabilization or regular homotopy. On the
other hand, many known examples of non-isotopic curves can be understood
in terms of braiding along Lagrangian annuli (or equivalently, Luttinger
surgery of the branched covers), leading to some intriguing open questions
about the topology of symplectic 4-manifolds versus that of K\"ahler surfaces.

If one prefers to adopt a more group theoretic point of view, it is possible
to use braid monodromy techniques to reformulate these questions in terms
of words in braid groups. For example, the classification of symplectic
4-manifolds reduces in principle to a (hard) question about factorizations
in the braid group, known as the Hurwitz problem.

In the following sections, we discuss these various questions and the
connections between them, starting from the point of view of symplectic
4-manifolds (in \S 2), then translating them in terms of plane branch
curves (in \S 3) and finally braid group factorizations (in \S 4).

\section{Topological questions about symplectic 4-manifolds}

\subsection{Classification of symplectic 4-manifolds}

Recall that a {\it symplectic manifold} is a smooth manifold equipped
with a 2-form $\omega$ such that $d\omega=0$ and $\omega\wedge\dots\wedge
\omega$ is a volume form. The first examples of compact symplectic manifolds
are compact oriented surfaces (taking $\omega$ to be an arbitrary area form),
and the complex projective space $\CP^n$ (equipped with the Fubini-Study
K\"ahler form). More generally, since any submanifold to which $\omega$
restricts non-degenerately inherits a symplectic structure, all complex
projective manifolds are symplectic. However, the symplectic category
is strictly larger than the complex projective category, as first
evidenced by Thurston in 1976 \cite{Th}. In 1994 Gompf used
the {\it symplectic sum} construction to prove that any finitely presented
group can be realized as the fundamental group of a compact symplectic
4-manifold \cite{Go1}.

An important problem in symplectic topology is to understand the hierarchy
formed by the three main classes of compact oriented 4-manifolds: (1)
complex projective, (2) symplectic, and (3) smooth. Each class is a proper
subset of the next one, and many obstructions and examples are known,
but we are still very far from understanding
what exactly causes a smooth 4-manifold to admit a symplectic structure, or
a symplectic 4-manifold to admit an integrable complex structure.

One of the main motivations to study symplectic 4-manifolds
is that they retain some (but not all) features of complex projective
manifolds: for example the structure of their Seiberg-Witten invariants,
which in both cases are non-zero and count embedded (pseudo)holomorphic
curves
\cite{Ta1,Ta2}. At the same time, every compact oriented smooth 4-manifold
with
$b_2^+\ge 1$ admits a ``near-symplectic'' structure, i.e.\ a closed 2-form
which vanishes along a union of circles and is symplectic over the
complement of its zero set \cite{GK,Ho1}; and it appears that some
structural properties of symplectic manifolds carry over to the world of
smooth 4-manifolds (see e.g.\ \cite{Ta3,Asinglp}).
\medskip

Although the question of determining which smooth 4-manifolds admit symplectic
structures and how many is definitely an essential one, it falls outside
of the scope of this paper. Rather, our goal will be to obtain information
on the richness of the symplectic category, especially when compared to the
complex projective category. 

We will restrict ourselves to the class of {\it integral} compact symplectic
4-manifolds, i.e.\ we will assume that the cohomology class $[\omega]\in
H^2(X,\R)$ is the image of an element of $H^2(X,\Z)$. This does not
place any additional restrictions on the diffeomorphism type of $X$, but
makes classification a discrete problem (by Moser's stability theorem,
deformations that keep $[\omega]$ constant are induced by ambient
isotopies). 

By integrating the Chern classes of the tangent bundle and the symplectic
class over the fundamental cycle $[X]$, one obtains various classical topological
invariants: the Chern numbers $c_1^2$ ($=2\chi+3\sigma$) and $c_2$ ($=\chi$),
the symplectic volume $[\omega]^2$, and $c_1\cdot[\omega]$. Hence, the first
question we will ask is:

\begin{question}\label{q:classif}
Can one classify all integral compact symplectic 4-manifolds with given
values of $(c_1^2,c_2,c_1\cdot[\omega],[\omega]^2)$ $($and
a given fundamental group$)$?
\end{question}

This question contains the {\it geography problem}, i.e.\ the question of
determining which Chern numbers can be realized by compact symplectic
4-manifolds. In some specific cases, Taubes' results on Seiberg-Witten
invariants seriously constrain the list of possibilities. For example
we have the following result \cite{Ta1}:

\begin{theorem}[Taubes]
Let $(X,\omega)$ be a compact symplectic 4-manifold with $b_2^+\ge 2$.
Then $c_1\cdot[\omega]\le 0$. Moreover, if $X$ is minimal
(i.e.\ does not contain an embedded symplectic sphere
of square $-1$), then $c_1^2\ge 0$.
\end{theorem}

A lot is also known about
the case $c_1^2=0$, from Seiberg-Witten theory
and from various surgery constructions. For example,
infinite families of simply connected symplectic 4-manifolds
homeomorphic but not diffeomorphic to elliptic surfaces have
been constructed (see e.g.\ \cite{Go1,FS1}).
However, when $c_1^2>0$ very little is known, and many important
questions remain open. For example it is unknown whether the
Bogomolov-Miyaoka-Yau inequality $c_1^2\le 3c_2$, which constrains the Chern
numbers of complex surfaces of general type, holds for symplectic
4-manifolds.

\subsection{Lefschetz fibrations and stabilization by symplectic sums}

One possible approach to the classification of symplectic 4-manifolds
is via {\it symplectic Lefschetz fibrations}, as suggested by Donaldson.
After blowing up a certain number of points, every compact
integral symplectic 4-manifold can be realized as the total space of a 
fibration over $S^2$ whose fibers are compact Riemann surfaces, finitely
many of which present a nodal singularity~\cite{Do2}. 
Conversely, the total space of such a Lefschetz fibration is
a symplectic 4-manifold \cite{GS}. If one could classify symplectic
Lefschetz fibrations, then an answer to Question
\ref{q:classif} would follow.

When the fiber genus is 0 or 1, the classification of Lefschetz fibrations
is a classical result; in particular, these fibrations are all
holomorphic \cite{Mo1}.
For genus~2, Siebert and Tian have proved holomorphicity
under assumptions of irreducibility of the singular fibers
and transitivity of the monodromy \cite{ST}, but in general there are
non-holomorphic examples \cite{OS}, and the complete classification is not
known. However, the situation simplifies if we ``stabilize'' by repeatedly
performing fiber sums with a specific holomorphic fibration $f_0$
(the fibration obtained by blowing up a pencil of curves of bidegree
$(2,3)$ in $\CP^1\times \CP^1$). Then we have the following result
\cite{Agenus2}:

\begin{theorem}
For any genus 2 symplectic Lefschetz fibration $f:X\to S^2$, there exists
an integer $n_0$ such that, for all $n\ge n_0$, $f\# nf_0$ is isomorphic
to a holomorphic fibration.
\end{theorem}

In fact, given two genus 2 symplectic Lefschetz fibrations $f,f'$ with
the same numbers of singular fibers of each type (irreducible, reducible
with genus 1 components, reducible with components of genus 0 and 2),
for all large $n$ the fiber sums $f\#nf_0$ and $f'\#nf_0$ are isomorphic
\cite{Agenus2}. More generally,
as a corollary of a recent result of Kharlamov and Kulikov about braid
monodromy factorizations \cite{KK}, a similar result holds for all Lefschetz
fibrations with monodromy contained in the hyperelliptic mapping class
group. This leads to the following questions relative to the classification
of symplectic 4-manifolds up to stabilization by fiber sums:

\begin{question}\label{q:stiso1}
Does every symplectic Lefschetz fibration become isomorphic to a holomorphic
fibration after repeatedly fiber summing with certain standard holomorphic
fibrations?
\end{question}

\begin{question}\label{q:stiso2}
Let $X_1,X_2$ be two integral compact symplectic 4-manifolds with the same
$(c_1^2,c_2,c_1\cdot[\omega],[\omega]^2)$. Do $X_1$ and $X_2$ become
symplectomorphic after repeatedly performing symplectic sums with the
same complex projective surfaces (chosen among a finite collection of model
surfaces)?
\end{question}

\subsection{Luttinger surgery}
Many of the constructions used to obtain interesting
examples of non-K\"ahler symplectic 4-manifolds, such as symplectic sum,
link surgery, and symplectic rational blowdown, rely on the idea of cutting
and pasting elementary building blocks. We focus here on the construction
known as {\it Luttinger surgery} \cite{Lu}, which has been comparatively
less studied but can be used to provide a unified description of numerous
examples of exotic symplectic 4-manifolds.

Given an embedded Lagrangian torus $T$ in a symplectic 4-manifold
$(X,\omega)$ and a homotopically non-trivial embedded loop $\gamma\subset T$,
Luttinger surgery is an operation that consists in
cutting out from $X$ a tubular neighborhood of $T$, foliated by parallel
Lagrangian tori, and gluing it back in such a way that the new meridian loop
differs from the old one by a twist along the loop $\gamma$ (while
longitudes are not affected), yielding a new symplectic manifold
$(\tilde{X},\tilde{\omega})$.

More precisely, identify a neighborhood of $T$ in $X$ with the neighborhood
$T^2\times D^2(r)$ of the zero section in
$(T^*T^2,dp_1\wedge dq_1+dp_2\wedge dq_2)$, in such a way
that $\gamma$ is identified with the first factor in $T^2=S^1\times S^1$.
Let $\theta$ be a smooth circle-valued function on the annulus
$A=D^2(r)\setminus D^2(\frac{r}{2})$ such that $\partial \theta/\partial
p_2=0$, and representing the generator of $H^1(A)=\Z$ (i.e., the value
of $\theta$ increases by $2\pi$ as one goes around the origin).
The diffeomorphism of $T^2\times A$ defined by 
$$\phi(q_1,q_2,p_1,p_2)=(q_1+\theta(p_1,p_2),q_2,p_1,p_2)$$
preserves the symplectic form, and so the manifold $$\tilde{X}=(X\setminus
T^2\times D^2(\tfrac{r}{2}))\cup_\phi (T^2\times D^2(r))$$
inherits a natural symplectic structure. For more details see \cite{Lu,ADK}.

By performing Luttinger surgery along suitably chosen Lagrangian tori,
one can e.g.\ transform a product $T^2\times \Sigma$ into any
surface bundle over $T^2$, or an untwisted fiber sum of Lefschetz
fibrations into a twisted fiber sum. Fintushel and Stern's symplectic
examples of knot surgery manifolds can also be obtained from complex
surfaces by Luttinger surgery. Although there is no good reason to
believe that the answer should be positive, the wide range of examples
which reduce to this construction makes it interesting to ask the
following question:

\begin{question}\label{q:luttinger}
Let $X_1,X_2$ be two integral compact symplectic 4-manifolds with the same
$(c_1^2,c_2,c_1\cdot[\omega],[\omega]^2)$. Is it always possible to obtain
$X_2$ from $X_1$ by a sequence of Luttinger surgeries?
\end{question}

In this question, as in Question \ref{q:stiso2} above, we do
not require the fundamental groups of $X_1$ and $X_2$ to be isomorphic.
This is because Luttinger surgery, like symplectic sum, can drastically
modify the fundamental group. Also, let us
mention that a positive answer to Question \ref{q:luttinger} essentially
implies a positive answer to Question \ref{q:stiso2}, as we shall see in
\S 4.

The symplectic sum construction can be used to build minimal simply connected
symplectic 4-manifolds with Chern numbers violating the Noether inequality,
and hence not diffeomorphic to any complex surface (see e.g.\
Theorem 10.2.14 in \cite{GS}). Many of these manifolds are homeomorphic
to (non-minimal) complex surfaces, but it is not clear at all whether it is
possible to obtain them by Luttinger surgeries. Given the very explicit
nature of the construction, these could be good test
examples for Question \ref{q:luttinger}.

\section{Isotopy questions about singular plane curves}

\subsection{Symplectic branched covers}

Let $X$ and $Y$ be compact oriented 4-manifolds, and assume that $Y$ carries
a symplectic form $\omega_Y$.

\begin{definition}
A smooth map $f:X\to Y$ is a {\em symplectic branched covering} if given any
point $p\in X$ there exist neighborhoods $U\ni p$, $V\ni f(p)$, and local
coordinate charts $\phi:U\to\C^2$
$($orientation-preserving$)$ and $\psi:V\to\C^2$
$($adapted to $\omega_Y$, i.e.\ such that $\omega_Y$ restricts positively
to any complex line in $\C^2)$, in which $f$ is given by one of:
\smallskip

$(i)$ $(x,y)\mapsto (x,y)$ $($local diffeomorphism$)$,

$(ii)$ $(x,y)\mapsto (x^2,y)$ $($simple branching$)$,

$(iii)$ $(x,y)\mapsto (x^3-xy,y)$ $($ordinary cusp$)$.
\end{definition}

These local models are the same as for the singularities of a generic 
holomorphic map from $\C^2$ to itself, except that the requirements on the
local coordinate charts have been substantially weakened.
The {\it ramification curve} $R=\{p\in X,\ \det(df)=0\}$ is a smooth
submanifold of $X$, and its image $D=f(R)$ is the {\it branch curve}, described
in the local models by the equations $z_1=0$ for $(x,y)\mapsto (x^2,y)$
and $27z_1^2=4z_2^3$ for $(x,y)\mapsto (x^3-xy,y)$. It follows from the
definition that $D$ is a singular
symplectic curve in $Y$.
Generically, its only singularities are
transverse double points, which may occur with either the complex
orientation or the opposite orientation, and complex cusps.
We have the following result \cite{Au2}:

\begin{proposition}\label{prop:au2}
Given a symplectic branched covering $f:X\to Y$, the manifold $X$ inherits
a natural symplectic structure $\omega_X$, canonical up to isotopy, in the
cohomology class $[\omega_X]=f^*[\omega_Y]$.
\end{proposition}

The symplectic form $\omega_X$ is constructed by adding to $f^*\omega_Y$
a small multiple of an exact form $\alpha$ with the property that, at
every point of $R$, the restriction of $\alpha$ to $\mathrm{Ker}(df)$ is
positive. Uniqueness up to isotopy follows from the
convexity of the space of such exact 2-forms and Moser's theorem.

Conversely, we can realize every integral compact symplectic 4-manifold
as a symplectic branched cover of $\CP^2$ \cite{Au2}:

\begin{theorem}\label{thm:au2}
Given an integral compact symplectic 4-manifold $(X^4,\omega)$ and an
integer $k\gg 0$, there exists a symplectic branched covering
$f_k:X\to\CP^2$, canonical up to isotopy if $k$ is sufficiently large.
\end{theorem}

The maps $f_k$ are built from suitably chosen triples of sections of
$L^{\otimes k}$, where $L\to X$ is a complex line bundle
such that $c_1(L)=[\omega]$. In the complex case, $L$ is an ample
line bundle, and a generic triple of holomorphic sections of
$L^{\otimes k}$ determines a $\CP^2$-valued map
$f_k:p\mapsto [s_0(p)\!:\!s_1(p)\!:\!s_2(p)]$.
In the symplectic case the idea is similar, but requires more analysis;
the proof relies on {\it asymptotically holomorphic} methods \cite{Au2}.

In any case, the natural symplectic structure induced on $X$ by the
Fubini-Study K\"ahler form and $f_k$ (as given
by Proposition \ref{prop:au2}) agrees with $\omega$ up to isotopy and
scaling (multiplication by~$k$).

Because for large $k$ the maps $f_k$ are canonical up to isotopy through
symplectic branched covers, the topology of $f_k$ and of its branch curve
$D_k$ can be used to define invariants of the symplectic
manifold $(X,\omega)$. Although the only generic singularities of the plane
curve $D_k$
are nodes (transverse double points) of either orientation and complex
cusps, in a generic one-parameter family of branched covers pairs of nodes
with opposite orientations may be cancelled or created. However, recalling
that a node of $D_k$ corresponds to the occurrence of two simple branch
points in a same fiber of $f_k$, the creation of a pair of nodes can only
occcur in a manner compatible with the branched covering structure, i.e.\
involving disjoint sheets of the covering.

It is worth mentioning that, to this date, there is no evidence suggesting
that negative nodes actually do occur in these high degree branch curves;
our inability to rule our their presence might well be a shortcoming of the
approximately holomorphic techniques, rather than an intrinsic feature of
symplectic 4-manifolds. So we will occasionally
consider the more conventional problem of understanding isotopy classes
of curves presenting only positive nodes and cusps, although most of the
discussion applies equally well to curves with negative nodes.
\medskip

Assuming that the topology of the branch curve is understood, 
the structure of $f$ is
determined by its {\it monodromy morphism} $\theta:\pi_1(\CP^2-D)\to S_N$,
where $N$ is the degree of the covering $f$. Fixing a base point $p_0\in
\CP^2-D$, the image by $\theta$ of a loop $\gamma$ in the complement of $D$
is the permutation of the fiber $f^{-1}(p_0)$ induced by the monodromy of $f$
along $\gamma$. (Since viewing this permutation as an element of $S_N$
depends on the choice of an identification
between $f^{-1}(p_0)$ and $\{1,\dots,N\}$, the morphism $\theta$ is only
well-defined up to conjugation by an element of $S_N$.) By Proposition
\ref{prop:au2}, the isotopy class of the branch curve $D$ and the monodromy
morphism $\theta$ determine completely the symplectic 4-manifold $(X,\omega)$
up to symplectomorphism.

The image by $\theta$ of a {\it geometric generator} of $\pi_1(\CP^2-D)$,
i.e.\ a loop $\gamma$ which bounds a small topological disc
intersecting $D$ transversely once, is a transposition (because
of the local model near a simple branch point). Since the image of $\theta$
is generated by transpositions and acts transitively on the fiber (assuming
$X$ to be connected), $\theta$ is a surjective group homomorphism. Moreover,
the smoothness of $X$ above the singular points of $D$ imposes certain
compatibility conditions on $\theta$. Therefore, not every singular plane
curve can be the branch curve of a smooth covering; in fact, the morphism
$\theta$, if it exists, is often unique (up to conjugation in $S_N$).
In the case of algebraic curves, this uniqueness property, which holds
except for a finite list of well-known counterexamples, is known as
Chisini's conjecture, and was essentially proved by Kulikov a few years
go \cite{Ku}.

The upshot of the above discussion is that, in order to understand symplectic
4-manifolds, it is in principle enough to understand singular plane curves.
Moreover, if the branch curve of a symplectic covering $f:X\to \CP^2$ happens
to be a complex curve, then the integrable complex structure of $\CP^2$ can be
lifted to an integrable complex structure on $X$, compatible with the
symplectic structure; this implies
that $X$ is a complex projective surface. So, considering the branched
coverings constructed in Theorem \ref{thm:au2}, we have:

\begin{corollary} \label{cor:au2}
For $k\gg 0$ the branch curve $D_k\subset\CP^2$ is isotopic to a complex
curve (up to node cancellations) if and only if $X$ is a complex projective
surface.
\end{corollary}

This motivates the study of the {\it symplectic isotopy problem} for
singular curves in $\CP^2$ (or more generally in other complex surfaces --
especially rational ruled surfaces, i.e.\ $\CP^1$-bundles over $\CP^1$).

\subsection{The symplectic isotopy problem} The symplectic isotopy problem
asks under which circumstances (assumptions on degree,
singularities, \dots) it is true that any symplectic curve is isotopic 
to a complex curve (by isotopy, we mean a continuous one-parameter family
of symplectic curves with the same singularities).
More generally, the goal is to understand isotopy
classes of symplectic curves with given singularities in a given homology
class. For example, considering only plane curves with positive nodes and cusps,
one may ask the following:

\begin{question}\label{q:curvclassif}
Given integers $(d,\nu,\kappa)$, can one classify all symplectic 
curves of degree $d$ in $\CP^2$ with $\nu$ nodes and $\kappa$ cusps,
up to symplectic isotopy?
\end{question}

If $D$ is the branch curve of an $N$-fold symplectic covering,
then the Chern classes of the symplectic manifold $(X,\omega)$ (with the
symplectic structure given by Proposition \ref{prop:au2}) are related to
the degree $d$ of $D$, its genus $g=\frac{1}{2}(d-1)(d-2)-\kappa-\nu$, and
its number of cusps via the formulas:
$$[\omega]^2=N,\quad c_1\cdot[\omega]=3N-d,\quad
c_1^2=g-1-\tfrac{9}{2}d+9N,\quad
c_2=2g-2+3N-\kappa.
$$
In particular, integrality constraints on the Euler-Poincar\'e
characteristic $\chi=c_2$ and signature $\sigma=\frac{1}{3}(c_1^2-2c_2)$
of $X$ imply that the degree $d$ must be even, and that the number of
cusps $\kappa$ must be a multiple of 3. The geography problem for symplectic
4-manifolds translates into a geography problem for symplectic branch
curves: for example, the Bogomolov-Miyaoka-Yau inequality $c_1^2\le 3c_2$
translates into the inequality $$\kappa\le \tfrac{5}{3}(g-1)+\tfrac{3}{2}d.$$
There are plane curves which violate this inequality, even in the algebraic
world: e.g.\ the branch curves of generic projections of irrational ruled
surfaces $\Sigma\times \CP^1$, where $\Sigma$ is a curve of genus $\ge 2$.
However, the open question is whether one can find branch curves which
violate this inequality and for which
the branched covering has $c_1^2\ge 0$. By the above remarks,
these cannot be isotopic to any complex curve.
\medskip

The symplectic isotopy problem is understood in various simple situations,
where it can be shown that every symplectic curve is isotopic to a complex
curve. The first results were obtained by Gromov \cite{Gr},
who used pseudoholomorphic
curves to prove that every smooth symplectic curve of degree $1$ or $2$ in
$\CP^2$ is isotopic to a complex curve. The idea of the argument is to equip
$\CP^2$ with an almost-complex structure $J=J_1$ such that the given
curve $C$ is $J$-holomorphic, and consider a smooth family of almost-complex
structures $(J_t)_{t\in [0,1]}$ interpolating between $J$ and the
standard complex structure $J_0$. By studying the deformation problem for
pseudoholomorphic curves, one can prove the existence of a smooth family of
$J_t$-holomorphic curves $C_t$ realizing an isotopy between $C=C_1$ and
an honest holomorphic curve $C_0$.

Successive improvements of this result have been obtained by Sikorav
(for smooth curves of degree $\le 3$), Shevchishin (degree $\le 6$), and
more recently Siebert and Tian \cite{ST}:

\begin{theorem}[Siebert-Tian]
Every smooth symplectic curve of degree $\le 17$ in $\CP^2$ is isotopic
to a complex curve.
\end{theorem}

\noindent
A similar result has also been obtained for smooth curves in
$\CP^1$-bundles over $\CP^1$ (assuming
$[C]\cdot[\mathrm{fiber}]\le 7$) \cite{ST}. It is expected that the isotopy
property remains true for smooth plane curves of arbitrarily large degree;
this would provide an answer to Question \ref{q:curvclassif} in the case
$\nu=\kappa=0$ (recall that all smooth complex curves of a given degree
are mutually isotopic).

The isotopy property is also known to hold in some simple cases for
curves with nodes and cusps in $\CP^2$ and $\CP^1$-bundles
over $\CP^1$,
as illustrated by the results obtained by Barraud, Shevchishin, and
Francisco. For example, we have the following results \cite{Sh,SF}:

\begin{theorem}[Shevchishin]
Any two irreducible nodal symplectic curves in $\CP^2$ of the same degree
and the same genus $g\le 4$ are symplectically isotopic.
\end{theorem}

\begin{theorem}[Francisco]
Let $C$ be an irreducible symplectic curve of degree $d$ and genus $0$
with $\kappa$ cusps and $\nu$ nodes in $\CP^2$, and assume that
$\kappa<d$.
Then $C$ is isotopic to a complex curve.
\end{theorem}

In general, we cannot expect the classification to be so simple, and there
are plenty of examples of symplectic curves which are not isotopic to any
complex curve. Perhaps the most widely known such examples are due to
Fintushel and Stern \cite{FS2}, who showed that elliptic surfaces contain
infinite families of
pairwise non-isotopic smooth symplectic curves representing a same homology
class. Similar results have also been obtained by Smith,
Etg\"u and Park, and Vidussi. However, if we consider singular curves with
cusp singularities, then these non-isotopy phenomena already arise in $\CP^2$.
In a non-explicit manner, it is clear that this must be the case, from
Corollary \ref{cor:au2}; however to this date the branch curves given by
Theorem \ref{thm:au2} for $k\gg 0$ have not been computed explicitly for
any non-complex examples. More explicitly, the following result is due to
Moishezon \cite{MoC} (see also \cite{ADK}):

\begin{theorem}[Moishezon]
For all $p\ge 2$, there exist infinitely many pairwise non-isotopic
singular symplectic curves of degree $d=9p(p-1)$ in $\CP^2$ with
$\kappa=27(p-1)(4p-5)$ cusps and $\nu=\frac{27}{2}(p-1)(p-2)(3p^2+3p-8)$
nodes, not isotopic to any complex curve.
\end{theorem}

\noindent
Moishezon's approach is purely algebraic (using braid monodromy
factorizations), and yields curves that are distinguished by the
fundamental groups of their complements \cite{MoC}. However a simpler
geometric description of his construction can be given in terms
of braiding constructions \cite{ADK}; cf.~\S 3.4.
\medskip

Questions \ref{q:classif} and \ref{q:curvclassif} are closely related to
each other, via Proposition~\ref{prop:au2} and Theorem \ref{thm:au2}.
Let us restrict ourselves to those plane curves which admit a compatible
symmetric group valued monodromy morphism, and assume that Chisini's
conjecture about the uniqueness of this morphism (excluding a specific
degree 6 curve) extends to the symplectic case. Then integral compact
symplectic 4-manifolds (up to
scaling of the symplectic form) are in one-to-one correspondence with
isotopy classes of singular symplectic plane branch curves up to an equivalence
relation which takes into account: (1)~the possibility of creating and
cancelling pairs of nodes, and (2) the dependence of the branch curve $D_k$
on the parameter $k$ in Theorem \ref{thm:au2}. This latter dependence, while
complicated and not quite understood in general, is nonetheless within
reach: see \cite{doubling} for a description of the relation between
$D_k$ and $D_{2k}$.

If one allows creations and cancellations of pairs of nodes, then the
classification problem becomes different, even considering only
curves with positive nodes and cusps. Indeed, it may happen that two
non-isotopic curves can be deformed into each other if one is allowed to
``push'' the curve through itself, creating or cancelling pairs of double
points in the process (such a deformation is called a {\it regular
homotopy}). In fact, in this case the classification becomes excessively
simple, as shown by the following result \cite{AKS}:

\begin{theorem}[A.-Kulikov-Shevchishin]\label{thm:reghom}
Any two irreducible symplectic curves with positive nodes
and cusps in $\CP^2$, of the same degree and with the same numbers of
nodes and cusps, are regular homotopic to each other.
\end{theorem}

What this means is that, when considering symplectic branch curves given
by Theorem \ref{thm:au2}, it is important to restrict oneself to {\it
admissible} regular homotopies, i.e.\ regular homotopies which are
compatible with the symmetric group valued monodromy morphism $\theta$.
When pushing the branch curve $D$ through itself, the two branches that
are made to intersect each give rise to a geometric generator of
$\pi_1(\CP^2-D)$. The requirement for admissibility of a node creation
operation is that the images by $\theta$ of these two geometric generators
should be transpositions acting on disjoint pairs of elements (i.e.,
the branching phenomena above the two intersecting branches of $D$ should
occur in different sheets of the covering). Thus the version of the isotopy
problem which naturally comes out of Theorem \ref{thm:au2} is the following:

\begin{question}\label{q:admreghom}
Given integers $(d,\nu,\kappa,N)$, can one classify all pairs $(D,\theta)$
where $D$ is a symplectic curve of degree $d$ in $\CP^2$ with $\nu_+$
positive nodes, $\nu_-$ negative nodes and $\kappa$ complex cusps,
$\nu_+-\nu_-=\nu$, and $\theta:\mbox{$\pi_1(\CP^2-D)$}\to S_N$ is a compatible
monodromy morphism, up to admissible regular homotopies?
\end{question}

\subsection{Hurwitz curves and stabilization}
In order to state the analogue of Question \ref{q:stiso1} for branch
curves, we need to introduce a slightly more restrictive category of curves,
known at {\it Hurwitz curves}. Roughly speaking, a Hurwitz curve in a ruled
surface is a curve which behaves like a generic complex curve
with respect to the ruling. In the case of $\CP^2$, we
consider the projection $\pi:\CP^2-\{(0:0:1)\}\to \CP^1$ given by
$(x:y:z)\mapsto (x:y)$, and we make the following definition:

\begin{definition}\label{def:hurwitz}
A curve $D\subset\CP^2$ $($not passing through $(0\!:\!0\!:\!1))$ is a
Hurwitz curve if $D$ is positively transverse to the fibers of $\pi$
at all but finitely many points, where $D$ is smooth and
non-degenerately tangent to the fibers.
\end{definition}

Hurwitz curves in $\CP^1$-bundles over $\CP^1$ can be defined similarly,
considering the projection to $\CP^1$ given by the bundle structure.

It is easy to see that any Hurwitz curve in $\CP^2$ can be made symplectic
by an isotopy through Hurwitz curves: namely, the image of any Hurwitz curve
by the rescaling map $(x:y:z)\mapsto (x:y:\lambda z)$ is a Hurwitz curve,
and symplectic for $|\lambda|\ll 1$. Moreover, Theorem \ref{thm:au2} can
be improved to ensure that the branch curves $D_k\subset\CP^2$ are Hurwitz
curves \cite{AK}. So, the discussion in \S\S 3.1--3.2 carries over to the
world of Hurwitz curves without modification.

After blowing up $\CP^2$ at $(0\!:\!0\!:\!1)$, we obtain the Hirzebruch
surface $\mathbb{F}_1$ (recall that 
$\mathbb{F}_n=\mathbb{P}(\mathcal{O}_{\mathbb{P}^1}\oplus
\mathcal{O}_{\mathbb{P}^1}(n))$), and any Hurwitz curve in $\CP^2$ becomes
a Hurwitz curve in $\mathbb{F}_1$, disjoint from the exceptional section.
The advantage of considering Hurwitz curves in Hirzebruch surfaces rather
than $\CP^2$ is that we can now introduce an operation of stabilization by
{\it pairwise fiber sum}. Namely, consider two Hurwitz curves $D_1\subset
\mathbb{F}_{n_1}$, $D_2\subset \mathbb{F}_{n_2}$, of the same degree $d$
relatively to the projection, i.e.\ such that
$[D_1]\cdot[F]=[D_2] \cdot [F]=d$, where $F$ is the fiber of the ruling.
Then, up to an isotopy among Hurwitz curves, we can assume that the
intersections of $D_1$ and $D_2$ with fixed fibers of the rulings coincide,
and we can smooth the normal crossing configuration $(\mathbb{F}_{n_1},D_1)
\cup_{\,\mathrm{fiber=fiber}}(\mathbb{F}_{n_2},D_2)$ into a
pair $(\mathbb{F}_{n},D)$, where $D$ is a Hurwitz curve in $\mathbb{F}_n$,
and $n=n_1+n_2$.

If a Hurwitz curve in $\mathbb{F}_n$ is a branch curve, then the ruling on
$\mathbb{F}_n$ lifts to a symplectic Lefschetz fibration on the branched
cover. Assuming that the symmetric group valued morphisms are compatible
(i.e.\ have the same restrictions to given fibers),
the fiber sum operation on the branch curves then corresponds to a fiber sum 
operation on the covers. Hence, the analogue of Questions \ref{q:stiso1} and
\ref{q:stiso2} asks
whether stabilization by fiber summing can be used to simplify the
classification of Hurwitz branch curves:

\begin{question}\label{q:stisohur}
Let $D_1,D_2$ be two Hurwitz curves in $\mathbb{F}_n$,
representing the same homology class and with the same numbers of
cusps and nodes. Assume that two compatible monodromy morphisms
$\theta_i:\pi_1(\mathbb{F}_n-D_i)\to S_N$ are given $(i\in\{1,2\})$,
and that there is a fiber $F\subset \mathbb{F}_n$ such that $F\cap D_1=
F\cap D_2$ and $\theta_{1|F-D_1}=\theta_{2|F-D_2}$. Is there a complex
curve $C\subset \mathbb{F}_{n'}$, compatible with the given monodromy
morphisms, such that the fiber sums $D_1\# C$ and $D_2\# C$ are isotopic
to each other as Hurwitz curves?
\end{question}

\noindent
To remain closer to the formulation of the questions in \S 2, one can
instead require the complex curve to be chosen among a finite list of standard
models (depending on the given monodromy morphisms $\theta_i$), but allow
several successive fiber sum operations. It is also interesting to ask
whether the final result of the fiber sum operations can always be assumed
to be isotopic to a complex curve.

Requiring compatibility with the given monodromy morphisms places
restrictions on the choice of the curve $C$, and makes the question more
difficult. Without this constraint the answer is known, and follows directly from a recent
result of Kharlamov and Kulikov about braid monodromy factorizations~\cite{KK}:

\begin{theorem}[Kharlamov-Kulikov]\label{thm:kk}
Let $D_1,D_2$ be two Hurwitz curves in $\mathbb{F}_n$,
representing the same homology class and with the same numbers of
cusps and nodes. Then there exists a smooth complex curve $C$ in
$\mathbb{F}_0=\CP^1\times\CP^1$ such that the fiber sums 
$D_1\# C$ and $D_2\# C$ are isotopic.
\end{theorem}

In this result, $C$ is in fact a smooth curve of bidegree $(a,b)$, where
$a=[F]\cdot[D_i]$, and $b\gg 0$ is chosen very large. Such a curve may be
obtained by smoothing a configuration consisting of $a$ sections
$\CP^1\times\{pt\}$ and $b$ fibers $\{pt\}\times \CP^1$. Hence, in this
case the fiber sum operation is equivalent to considering the union of $D_i$
with $b$ fibers of the ruling, and smoothing the intersections; a more
geometric formulation of Theorem \ref{thm:kk} is therefore:

\begin{theorem}
Let $D_1,D_2$ be two Hurwitz curves in $\mathbb{F}_n$,
representing the same homology class and with the same numbers of
cusps and nodes. Let $D'_i$ $(i\in\{1,2\})$ be the curve obtained
by adding to $D_i$ a union of $b$ generic fibers of the ruling,
intersecting $D_i$ transversely at smooth points, and smoothing out
all the resulting intersections. Then for all
large enough values of $b$ the Hurwitz curves $D'_1$ and $D'_2$
are isotopic.
\end{theorem}

This construction gives an answer to Question \ref{q:stisohur} in
the case of smooth curves (and coverings of degree $N=2$); it is
unclear whether the argument in \cite{KK} can be modified to
produce complex curves compatible with branched coverings of
degree $N\ge 3$.

\subsection{Braiding along Lagrangian annuli}
Let $D$ be a symplectic curve in a symplectic 4-manifold $Y$ (e.g.\
$Y=\CP^2$), possibly with singularities. It is often the
case that we can find an embedded Lagrangian annulus \hbox{$A\subset Y
\setminus D$,} with boundary contained in the smooth part of $D$. 
(This happens for example when a portion of $D$ consists of two cylinders
which run parallel to each other; then we can find a Lagrangian annulus
joining them).

In this situation, one can twist the curve $D$ along the
annulus $A$, to obtain a new symplectic curve $\tilde{D}$ which coincides with $D$ away
from $A$ \cite{ADK}. Namely, we can identify a neighborhood of $A$ with the product
$S^1\times (-1,1) \times D^2$, in such a way that $A=S^1\times\{0\}\times
[-\frac{1}{2},\frac{1}{2}]$ and a neighborhood of $\partial A$ in $D$
is $S^1\times (-1,1)\times \{\pm \frac{1}{2}\}$. (If we deform $D$ suitably,
then we may assume that the symplectic structure is the product one, but
this is not necessary). Then the curve $\tilde{D}$ is obtained from $D$ by
replacing $S^1\times (-1,1)\times \{\pm \frac{1}{2}\}$ by
$S^1\times\tilde\Gamma$, where $\tilde\Gamma=
\{(t,\,\pm \frac{1}{2}\exp(i\pi\chi(t))),\ t\in(-1,1)\}\subset (-1,1)\times
D^2$ and $\chi$ is a smooth function which equals $0$ near $-1$ and $1$
near $1$. This construction is called ``braiding'' because, forgetting
the $S^1$ factor, it replaces the trivial braid $(-1,1)\times \{\pm
\frac{1}{2}\}$ with the half-twist $\tilde\Gamma$.

Assume now that $D$ is the branch curve of an $N$-fold symplectic
covering $f:X\to Y$. Assume moreover that $f$ is ramified in the
same manner above the two boundary components of $A$, i.e.\ that
two of the $N$ lifts of $A$ have boundary contained in the
ramification curve $R$; then these two lifts together form an embedded
Lagrangian torus $T\subset X$, and we have the following result \cite{ADK}:

\begin{proposition}[A.-Donaldson-Katzarkov]
The symplectic 4-manifold $\tilde{X}$
obtained from $X$ by Luttinger surgery along the torus $T$ is the total
space of a natural symplectic branched covering $\tilde{f}:\tilde{X}\to Y$,
whose branch curve $\tilde{D}$ is the curve obtained from $D$ by
braiding along the annulus $A$.
\end{proposition}

Hence, the natural analogue of Question \ref{q:luttinger}
for singular plane curves is:

\begin{question}\label{q:braiding}
Let $D_1,D_2$ be two symplectic curves with positive nodes
and cusps in $\CP^2$, of the same degree and with the same numbers of
nodes and cusps. 
Is it always possible to obtain $D_2$ from $D_1$ by a sequence
of braiding operations along Lagrangian annuli?
\end{question}

As before, there is no good reason to believe that the answer should be
positive, except that most known examples of non-isotopic symplectic
curves seem to reduce to this construction. This is e.g.\ the case for
the Fintushel-Stern examples of non-isotopic smooth symplectic curves in
elliptic surfaces \cite{FS2}, which are obtained by braiding a disconnected
union of elliptic fibers, and for Moishezon's examples
of singular plane curves \cite{MoC,ADK}, which are obtained by braiding the
branch curve of the projection of a complex surface of general type.

\section{Questions about braid monodromy factorizations}

\subsection{The braid monodromy of a plane curve}

One of the main tools to study algebraic plane curves is the notion of
{\it braid monodromy}, which has been used extensively
by Moishezon and Teicher (among others) since the early 1980s in order to study the
branch curves of generic projections of complex projective surfaces
(see \cite{Te1} for a detailed overview).
Braid monodromy techniques apply equally well to the more general case of
Hurwitz curves in $\CP^2$ or more generally in rational ruled surfaces
(see Definition~\ref{def:hurwitz}).

Given a Hurwitz curve $D$ in $\CP^2$,
the projection $\pi:\CP^2-\{(0:0:1)\}\to\CP^1$ 
makes $D$ a singular branched cover of $\CP^1$,
of degree $d=\deg D$. Each fiber of $\pi$ is a
complex line $\ell\simeq \C\subset\CP^2$, and if $\ell$ does not pass
through any of the singular points of $D$ nor any of its vertical
tangencies, then $\ell\cap D$ consists of $d$ distinct points.
We can trivialize the fibration $\pi$ over an affine subset
$\C\subset\CP^1$, and define
the {\it braid monodromy morphism}
$$\rho:\pi_1(\C-\mathrm{crit}(\pi_{|D}))\to B_d.$$
Here $B_d$ is the Artin braid group on $d$ strings (the fundamental group
of the configuration space $\mathcal{X}_d$ 
of $d$ distinct points in $\C$), and for any loop
$\gamma$ the braid $\rho(\gamma)$ describes the
motion of the $d$ points of $\ell\cap D$ inside the
fibers of $\pi$ as one moves along the loop $\gamma$.

Equivalently, choosing an ordered system of arcs generating the free group
$\pi_1(\C-\mathrm{crit}(\pi_{|D}))$, one can express the braid monodromy
of $D$ by a {\it factorization} $$\Delta^2=\prod_{i} \rho_i$$ of the central
element $\Delta^2$ (representing a full rotation by $2\pi$) in $B_d$, where
each factor $\rho_i$ is the monodromy
around one of the special points (cusps, nodes, tangencies) of $D$.

The monodromy around a tangency point is a {\it half-twist} exchanging two
strands, i.e.\ an element conjugated to one of the standard generators of
$B_d$; the monodromy around a positive (resp.\ negative) node is the square
(resp.\ the inverse of the square) of a half-twist; and the monodromy around
a cusp is the cube of a half-twist. Hence, we are interested in
factorizations of $\Delta^2$ into products of powers of half-twists.

A same Hurwitz curve can be described by different factorizations of
$\Delta^2$ in $B_d$: switching to a different ordered system of generators
of $\pi_1(\C-\mathrm{crit}(\pi_{|D}))$ affects the collection of factors
$\langle \rho_1,\dots,\rho_r\rangle $ by a sequence of {\it Hurwitz moves},
i.e.\ operations of the form
$$\langle \rho_1,\,\cdots,\rho_i,\rho_{i+1},\,\cdots,\rho_r\rangle \,
\longleftrightarrow\, \langle \rho_1,\,\cdots,(\rho_i\rho_{i+1}\rho_i^{-1}),
\rho_i,\,\cdots,\rho_r\rangle;
$$
and changing the trivialization of the reference fiber
$(\ell,\ell\cap D)$ of $\pi$ (i.e.\ its identification with
the base point in $\mathcal{X}_d$)
affects braid monodromy by a {\it global conjugation}
$$\langle\rho_1,\,\cdots,\rho_r\rangle \,\longleftrightarrow\,
\langle b^{-1}\rho_1 b,\,\cdots,b^{-1}\rho_r b\rangle.
$$
For Hurwitz curves whose only singularities are cusps and nodes (of either
orientation), the braid monodromy factorization determines the isotopy
type completely (see for example \cite{KK}). Hence,
determining whether two given Hurwitz curves are isotopic 
is equivalent to determining whether two given factorizations of
$\Delta^2$ coincide up to Hurwitz moves and global conjugation.
In this language the isotopy problem for Hurwitz curves in $\CP^2$
becomes:

\begin{question}
Given integers $(d,\nu,\kappa)$, can one classify all factorizations of the
central element $\Delta^2$ of $B_d$ into a product of 
$\tau=d(d-1)-2\nu-3\kappa$ half-twists, $\nu$
squares of half-twists, and $\kappa$ cubes of half-twists, up to Hurwitz
moves and global conjugation?
\end{question}

If our goal is to consider only branch curves of
symplectic coverings (rather than arbritrary plane Hurwitz
curves), then we need
to look specifically for factorizations in which the factors belong to
the {\it liftable braid group}, i.e.\ the subgroup of $B_d$ consisting of
all braids compatible with given branched covering data.

More precisely, assume that a Hurwitz curve $D$ is the branch curve of
a symplectic branched covering $f:X\to\CP^2$.
The fibers of $\pi$ form a pencil of lines on $\CP^2$, 
whose preimages by $f$ equip $X$ with a structure of symplectic
Lefschetz pencil.
By restricting the monodromy of the covering to
a fiber $\ell$ of $\pi$, we obtain a symmetric group valued morphism
$$\theta:\pi_1(\ell-(\ell\cap D))\to S_N,$$
which describes how to realize a fiber of the Lefschetz pencil as a branched
covering of a fiber of $\pi$. The braid group acts on $\pi_1(\ell-(\ell\cap
D))$ by automorphisms; call $b_*$ the automorphism induced by the braid $b$.
Then the liftable braid group is
$$LB_d(\theta)=\{b\in B_d,\ \theta\circ b_*=\theta\}.$$
Equivalently, recall that $B_d$ is the fundamental group of the space
$\mathcal{X}_d$ of configurations of $d$ distinct points in $\C$, 
and consider the configuration space $\tilde{\mathcal{X}}_d$ whose elements
are pairs $(\Pi, \sigma)$, where $\Pi$ is a set of $d$ distinct points in
$\C$, and $\sigma$ is a surjective group homomorphism from $\pi_1(\C-\Pi)$
to $S_N$ mapping generators to transpositions.
The projection $(\Pi, \sigma)\mapsto \Pi$ is a finite covering,
and taking $\tilde{*}=(\ell\cap D,\theta)$ as base point
we have $LB_d(\theta)=\pi_1(\tilde{\mathcal{X}}_d,\tilde{*})$.

\subsection{Stabilization and partial conjugation}

The main feature which makes braid groups algorithmically manageable is
the {\it Garside property}. Namely, if we consider the semigroup of {\it
positive braids} $B_d^+$, defined by the same generators $(\sigma_i)_{1\le i\le
d-1}$ and relations ($\sigma_i\sigma_{i+1}\sigma_i=\sigma_{i+1}\sigma_i
\sigma_{i+1}$ $\forall i$ and $\sigma_i\sigma_j=\sigma_j\sigma_i$ $\forall
|i-j|\ge 2$) as $B_d$ but without allowing inverses of the generators, then
we have the following property \cite{Gar}:

\begin{theorem}[Garside]
The natural homomorphism $i:B_d^+\to B_d$ is an embedding.
\end{theorem}

In other terms, if two positive words in the generators of the braid group
represent the same braid, then they can be transformed into each other by
repeatedly using the defining relations, without ever introducing
the inverses of the generators. Garside's other fundamental observation is
that for any $b\in B_d$ there exists an integer $k$ and a positive braid
$\beta\in B_d^+$ such that $\Delta^{2k}b=i(\beta)$ \cite{Gar}. These
properties make it possible to obtain solutions to the word and conjugacy
problems (see also \cite{BKL} for a more modern approach); they also yield
a stable classification of braid group factorizations \cite{KK}.

Namely, let $\mathcal{F}_0$ be the standard factorization
$\Delta^2=(\sigma_1\cdot\ldots\cdot\sigma_{d-1})^d$ in $B_d$, and  
say that two factorizations $\mathcal{F}=(\rho_1\cdot\ldots\cdot \rho_r)$, 
$\mathcal{F}'=(\rho'_1\cdot
\ldots\cdot \rho'_r)$ have the same number of factors of each type if
they have the same number of factors $r$ and there is
a permutation $\sigma\in S_r$ such that $\rho_i$ is conjugated to
$\rho'_{\sigma(i)}$ for all $i=1,\dots,r$. 
Then the following result holds \cite{KK}:

\begin{theorem}[Kharlamov-Kulikov]\label{thm:kk2}
Let $\mathcal{F}$ and $\mathcal{F}'$ be two factorizations of the same
element in $B_d$, with the same numbers of factors of each type.
Then there exists an integer
$n$ such that the factorizations $\mathcal{F}\cdot (\mathcal{F}_0)^n$ and
$\mathcal{F}'\cdot (\mathcal{F}_0)^n$ are equivalent under Hurwitz moves.
\end{theorem}

\noindent (Here $\mathcal{F}\cdot(\mathcal{F}_0)^n$ is the factorization
consisting of the factors of $\mathcal{F}$ followed by those of
$\mathcal{F}_0$ repeated $n$ times).

Theorem \ref{thm:kk} follows from this result by specifically considering
factorizations of $\Delta^{2n}$ whose factors are powers of half-twists
and observing that $\mathcal{F}_0$ is the braid monodromy factorization of
a smooth algebraic plane curve.
However, considering that the factors in $\mathcal{F}_0$ generate the
entire braid group, of which $LB_d(\theta)$ is a proper subgroup as soon
as the degree $N$ of the covering is at least 3, one is prompted to ask
the following question:

\begin{question} \label{q:stfact}
Given a symmetric group valued morphism $\theta$, does a statement
similar to Theorem \ref{thm:kk2} hold for factorizations in $LB_d(\theta)$?
\end{question}

Assuming that the factorization in $LB_d(\theta)$ playing the role of
the standard factorization $\mathcal{F}_0$ in this statement can be
realized as the braid monodromy of an algebraic curve,
a positive answer to this question would imply positive answers
to Questions \ref{q:stiso2} and \ref{q:stisohur}.
\medskip

Finally, the last question we will consider is that of {\it partial
conjugation} of braid factorizations. Namely, given a factorization
$\mathcal{F}$ with factors $\rho_1,\dots,\rho_r$, integers
$1\le p< q\le r$, and a braid $b$ such that 
$\prod_{p\le i\le q} \rho_i$ commutes with $b$, we can form a new
factorization $\mathcal{F}'$, with factors
$\rho_1,\dots,\rho_{p-1},(b^{-1}\rho_p b),\dots,$ $(b^{-1}\rho_q b),\rho_{q+1},
\dots,\rho_r$: the {\it partial conjugate} of $\mathcal{F}$ by $b$.

\begin{lemma}\label{l:pconj}
If the element $b$ belongs to the subgroup of $B_d$ generated by
$\rho_p,\dots,\rho_q$, and if $\prod_{p\le i\le q} \rho_i$ is central
in this subgroup, then $\mathcal{F}'$ is equivalent to
$\mathcal{F}$ under Hurwitz moves.
\end{lemma}

The proof is easy, and relies on the same trick as in Lemma 6
of~\cite{Agenus2}. On the other hand, if $b$ does not
belong to the subgroup generated by the factors of $\mathcal{F}$,
then we can get interesting examples of
inequivalent factorizations; this is e.g.\ how Moishezon's examples
\cite{MoC} are constructed. 

\begin{question}\label{q:pconj}
Are any two factorizations of the same
element in $B_d$ $($resp.\ $LB_d(\theta)$$)$, with the
same numbers of factors of each type,
equivalent under Hurwitz moves and partial conjugations by elements of
$B_d$ $($resp.\ $LB_d(\theta)$$)$?
\end{question}

A positive answer to this question (for factorizations in $LB_d(\theta)$) 
would imply that Questions \ref{q:luttinger}
and \ref{q:braiding} also admit positive answers.
In fact, if one specifically
considers factorizations of $\Delta^2$ into a product
of powers of half-twists in $LB_d(\theta)$, then Questions
\ref{q:luttinger} and \ref{q:pconj} are almost equivalent. This is
because, given an arbitrary Lagrangian torus $T$ in a symplectic 4-manifold,
one can build a symplectic Lefschetz pencil for which $T$ fibers above
an embedded loop $\delta$ in $\CP^1$ and intersects each fiber above
$\delta$ in a simple closed loop $\gamma$. Luttinger surgery
along $T$ then amounts to a partial conjugation of the monodromy of the
pencil by the Dehn twist about $\gamma$, and considering branched coverings
of $\CP^2$ instead of Lefschetz pencils it should also amount to a partial
conjugation of the braid monodromy of the branch curve.

Moreover, a positive answer to Question \ref{q:pconj} also implies a
positive answer to Question \ref{q:stfact} (and hence to Questions
\ref{q:stiso2} and \ref{q:stisohur}), at least provided that there
exists an algebraic plane branch curve whose braid monodromy
generates the liftable braid subgroup $LB_d(\theta)$. The existence
of such a factorization $\mathcal{F}_{0,\theta}$ is rather likely, and
examples should be relatively easy to find, although the question has
not been studied. Assuming that this is the case, given two factorizations
$\mathcal{F}_1,\mathcal{F}_2$ in $LB_d(\theta)$ with the same numbers of
factors of each type, the factors in $\mathcal{F}_1\cdot
\mathcal{F}_{0,\theta}$ and $\mathcal{F}_2\cdot \mathcal{F}_{0,\theta}$
generate $LB_d(\theta)$, and hence, by Lemma \ref{l:pconj},
any partial conjugation operation performed on $\mathcal{F}_1\cdot
\mathcal{F}_{0,\theta}$ is equivalent to a sequence of Hurwitz moves.
So, if $\mathcal{F}_1\cdot
\mathcal{F}_{0,\theta}$ and $\mathcal{F}_2\cdot \mathcal{F}_{0,\theta}$
are equivalent under Hurwitz moves and partial conjugations then they are
equivalent under Hurwitz moves only.

\bigskip

\noindent {\it Note added in proof:}
Questions 2.4 and 2.5 have now essentially
been solved. The reader is referred to: D. Auroux, {\sl A stable
classification of Lefschetz fibrations}, Geom. Topol. {\bf 9}
(2005), 203--217.

\end{document}